\documentclass[a4paper, 12pt]{article}
\usepackage[cp1251]{inputenc}
\usepackage[russian]{babel}
\usepackage{amsmath, amssymb, amsfonts, amsthm}
\usepackage{indentfirst}
\usepackage[pdftex]{graphicx}
\usepackage[left=3cm, right=2cm, vmargin=2cm]{geometry}
\usepackage{setspace}
\usepackage[framemethod=1]{mdframed}
\usepackage{hyperref}
\usepackage{algorithm}
\usepackage{algorithmic}
\floatname{algorithm}{\normalfont\textrm{А л г о р и т м}}
\DeclareMathOperator*{\argmin}{arg\,min}

\graphicspath{{images/}}
\newtheorem{theorem}{Теорема}[section]

\newtheorem{propose}[theorem]{Предложение}
\newtheorem{example}{Пример}
\newtheorem{corollary}[theorem]{Следствие}
\newtheorem{remark}[theorem]{Замечание}
\begin{document}
\begin{center}
{\bf АДАПТИВНЫЕ МЕТОДЫ ГРАДИЕНТНОГО ТИПА ДЛЯ ЗАДАЧ ОПТИМИЗАЦИИ С ОТНОСИТЕЛЬНОЙ ТОЧНОСТЬЮ И ОСТРЫМ МИНИМУМОМ}
\end{center}
\begin{center}
Ф.\,С.\,Стонякин, С.\,С.\,Аблаев, И.\,В.\,Баран
\end{center}

\section{Введение}

Предлагаемая работа посвящена исследованию некоторых методов градиентного типа для задач выпуклой минимизации с относительной точностью, а также субградиентных методов для задач с острым минимумом. Первый рассматриваемый в настоящей работе класс задач имеет вид
\begin{equation}\label{eq1}
f(x) \rightarrow \min_{x \in Q},
\end{equation}
где $Q$~--- выпуклое замкнутое подмножество в $\mathbb{R}^n$, $f$--- выпуклая положительно однородная функция степени 1. Такая постановка задачи возникла в работах Ю.\,Е.\,Нестерова \cite{Nest_diss, NesRound, NesUnco}. Будем говорить, что значение $f(\widehat{x})$, вычисленное в точке $\widehat{x}$ допустимого множества задачи $Q$, приближает минимальное значение $f^{\ast} > 0$ с относительной точностью $\gamma > 0$, если
$$f(\widehat{x}) \leq (1 + \gamma) f^{\ast} .$$
Отметим, что в работах \cite{Nest_diss, NesRound, NesUnco} обоснована целесообразность рассмотрения использования относительной точности для понимания качества приближённого решения некоторых прикладных задач (например, проектирование механических конструкций), а также адекватность предположения положительной однородности $f$. Вслед за упомянутыми выше работами Ю.\,Е.\,Нестерова (см., например, главу 6 из~\cite{Nest_diss}) для задач с относительной точностью \eqref{eq1} сделаем следующее предположение:
\begin{equation}\label{eq2}
\text{dom} \; f \equiv \mathbb{R}^n, \quad 0 \in \text{int} \; \partial f(0).
\end{equation}
Иными словами, $f$~--- опорная функция некоторого выпуклого компактного множества, у которого 0 --- внутренняя точка. Как показано в \cite{Nest_diss, NesRound, NesUnco}, такое предположение естественно и оно покрывает широкий класс задач, хотя в общем случае и нет возможности гарантировать, что $f$ удовлетворяет \eqref{eq2}.

В разделе 2 разработан аналог ускоренного универсального градиентного метода для задач выпуклой положительно однородной минимизации и получена оценка количества итераций \eqref{eq-complexity-alg2}, после которых гарантированно выполняется предложенный адаптивный критерий остановки. Выполнение этого критерия остановки заведомо обеспечивает достижение требуемой относительной точности точки выхода алгоритма по целевому функционалу. Доказанная оценка \eqref{eq-complexity-alg2} сублинейна, но за счёт адаптивного критерия остановки потенциально возможно наблюдать более высокую скорость сходимости, что показано с помощью экспериментов в примере \ref{example1}.

Однако естественный интерес представляет вопрос о возможности гарантировать линейную скорость сходимости для негладких задач, пусть и с дополнительными предположениями. Поэтому далее в разделах 3 и 4 были исследованы субградиентные методы для негладких задач с теоретическими результатами о линейной скорости сходимости. Хорошо известно, что это возможно при дополнительных предположениях, что минимум острый и доступно оптимальное значение $f^{*}$. В работе \cite{Pol_1969} при таких условиях предложен способ выбора шага субградиентного метода, который гарантирует линейную скорость сходимости по аргументу. Отметим, что условие острого минимума верно, например, для задачи проектирования точки на выпуклый компакт или задачи о нахождении общей точки системы множеств. В настоящей работе предлагаются варианты субградиентных методов, один из которых гарантирует линейную скорость сходимости на классе слабо $\beta$-квазивыпуклых липшицевых функционалов ($\beta \in (0; 1]$), а другой --- на классе обычных квазивыпуклых гёльдеровых функционалов с острым минимумом. Напомним, что слабая $\beta$-квазивыпуклость введена несколько лет назад в \cite{Hardt}. Это понятие обобщает условие выпуклости. Например, невыпуклая функция одной переменной $f(x) = |x|(1-e^{-|x|})$ слабо $1$-квазивыпукла.

Опишем кратко структуру и новые результаты (вклад) предлагаемой статьи. Работа состоит из введения (раздел 1), заключения (раздел 5) и трёх основных разделов. В разделе 2 рассмотрен класс задач минимизации выпуклых положительно однородных функций с гёльдеровым градиентом (субградиентом при $\nu = 0$). Предложен специальный подход к выбору параметров и правилу остановки адаптивного метода подобных треугольников на классе задач минимизации выпуклых положительно однородных функций. Это позволило обосновать применимость соответствующего аналога универсального градиентного метода для задач с относительной точностью и доказать оптимальную оценку его скорости сходимости на выделенном классе задач (теорема \ref{then3}). Приведён пример результатов численных экспериментов (пример \ref{example1}), иллюстрирующих возможность повышения скорости сходимости за счёт предлагаемого адаптивного правила остановки \eqref{krost}. При этом сублинейная теоретическая оценка скорости сходимости алгоритма \ref{alg2} побудила нас исследовать также и субградиентные методы в случаях, когда  можно гарантировать линейную скорость сходимости. В разделе 3 предложен вариант субградиентного метода для задач минимизации слабо $\beta$-квазивыпуклых функционалов в случае острого минимума и доказан результат о линейной скорости сходимости для липшицевых функций указанного типа (теорема \ref{thmbetaquasiconv} и следствие \ref{corol}). В  разделе 4 предложен вариант субградиентного метода, который сходится с линейной скоростью на классе задач минимизации квазивыпуклых гёльдеровых функционалов с острым минимумом. Рассмотрено обобщение этого результата на класс задач со слабым острым минимумом (пример \ref{example3}). Сформулированы следствия из теорем \ref{thmbetaquasiconv} и \ref{th5} для задач с относительной точностью на классе выпуклых положительно однородных функционалов степени 1 (следствия \ref{corol2} и \ref{corol3}).

Обозначим через $\|\cdot \|_2$ стандартную евклидову норму в пространстве $\mathbb{R}^n$. Напомним вспомогательный результат из \cite{Nest_diss} (теорема 6.1.1).
\begin{theorem}\label{then1}
Пусть для некоторого $\gamma_0>0$ верно
\begin{equation}\label{eqineq}
f(x) \geq \gamma_0 \| x \|_2 \quad \forall\,x \in Q.
\end{equation}
Тогда для всякого точного решения $x_{*}$ задачи \eqref{eq1} и $ f^{*}=f(x_*)$ верно неравенство
$$\|x_0 - x_{*} \|_2 \leq \frac{2}{\gamma_0} f^{*} \leq \frac{2}{\gamma_0} f(x_0).$$
\end{theorem}
Заметим, что предположение \eqref{eqineq} теоремы \ref{then1} есть следствие условия \eqref{eq2}. Для результатов (теорема \ref{then3}, следствия \ref{corol2} и \ref{corol3}) о задачах с относительной точностью по аналогии с \cite{Nest_diss} условимся в качестве начальной точки $x_0$ использовать евклидову проекцию начала координат на допустимое множество задачи $Q$: $\|x_0\|_2 :=\min\limits_{x \in Q} \{ \| x \|_2 \}$.

\section{Универсальный метод для задач минимизации выпуклого положительно однородного функционала с относительной точностью}

Предположим, что для градиента целевого функционала $f$ задачи \eqref{eq1} выполняется  условие Гёльдера, т.е. существует $\nu \in [0, 1]$, такое, что при некотором $L_{\nu} \in (0; + \infty)$ верно неравенство
\begin{equation}\label{eqHolderGrad}
\| \nabla f(x) - \nabla f(y) \|_2 \leq L_{\nu} \| x -y \|_2^{\nu} \quad \forall\,x,y\in Q.
\end{equation}
Если $\nu = 0$, то под $\nabla f(x)$ в \eqref{eqHolderGrad} можно понимать произвольный субградиент $f$ в точке $x$. При указанном допущении  \eqref{eqHolderGrad}, как известно \cite{Gasnikov_book, 4}, для произвольных $x, y \in Q$ верно неравенство
\begin{equation}\label{eqgensmooth}
0 \leq f(x) - f(y) - \langle \nabla f(y), x - y \rangle \leq \frac{L(\delta)}{2} \|x - y \|_2 + \delta
\end{equation}
при некотором фиксированном $\delta > 0$ и $$ L(\delta) = L_{\nu} \left [ \frac{L_{\nu}}{2 \delta} \frac{1-\nu}{1+\nu}\right ]^{\frac{1-\nu}{1+\nu}}.$$

Опишем универсальный метод, который работает на выделенном классе задач \eqref{eq1} с реализацией в ходе работы адаптивной настройки оценки качества выдаваемого решения с относительной точностью на параметр гладкости $L(\delta)$. Идея универсальных методов состоит в том, что используемый в оценке качества решения потенциально большой параметр гладкости $L = L(\delta)$ заменяется адаптивно подбираемым на итерациях множеством значений его локальных аналогов $L_{k+1}$ ($k \geq 0$), что может повышать качество выдаваемого решения по сравнению с теоретической оценкой, содержащий глобальное значение $L$ \cite{4}. Мы отправляемся от следующего адаптивного варианта метода подобных треугольников~\cite{3}. Сформулируем его вариант для задачи минимизации выпуклой функции $f$ с условием обобщённой гладкости \eqref{eqgensmooth}.

\begin{algorithm}[ht]
\caption{Адаптивный вариант метода подобных треугольников для задач с условием обобщённой гладкости.}
\label{alg2}
\begin{algorithmic}[1]
\REQUIRE $x_0$~--- начальная точка, $N$~--- количество шагов, $\{\delta_{k+1}\}_{k=0}^{N-1}$~--- последовательность положительных чисел и $L_0 > 0$;
\STATE \textbf{0-шаг:}
$
y_0 := x_0,\,
u_0 := x_0,\,
L_1 := \frac{L_0}{2},\,
\alpha_0 := 0,\,
A_0 := \alpha_0;
$
\FOR{$k = 1, ...$}
\STATE Найти наибольший корень $\alpha_{k+1}$ квадратного уравнения
$$A_k + \alpha_{k+1} = L_{k+1}\alpha^2_{k+1}, \quad A_{k+1} := A_k + \alpha_{k+1}, \quad y_{k+1} := \frac{\alpha_{k+1}u_k + A_k x_k}{A_{k+1}};$$
$$u_{k+1} := { \argmin\limits_{x \in Q}} \; \left \{ \frac{1}{2} \| x - u_k \|_2^2 + \alpha_{k+1} \left \langle \nabla f(y_{k+1}), x - y_{k+1} \right \rangle \right \};$$
$$x_{k+1} := \frac{\alpha_{k+1}u_{k+1} + A_k x_k}{A_{k+1}};$$
\IF{$f(x_{k+1}) \leq f(y_{k+1}) + \left \langle \nabla f (y_{k+1}), x_{k+1} - y_{k+1} \right \rangle + \frac{L_{k+1}}{2}\|x_{k+1} - y_{k+1}\|_2^2 + \delta_{k+1}$}
\STATE $L_{k+2} := \frac{L_{k+1}}{2}$ и перейти к следующему шагу.
\ELSE
\STATE $L_{k+1} := 2L_{k+1}$ и повторить текущий шаг.
\ENDIF
\ENDFOR
\end{algorithmic}
\end{algorithm}

Для алгоритма \ref{alg2} имеет место следующая (см. \cite{3}, теорема 2)
\begin{theorem}\label{then2}
Пусть $\|x_{*} - x_0 \|_2 \leq R\sqrt{2}$ для некоторой постоянной $R > 0$, где $x_0$~--- начальная точка, а $x_*$~--- ближайшая точка минимума $f$ к $x_0$ по евклидовой норме. Тогда для алгоритма~\ref{alg2} выполнено следующее неравенство:
\begin{equation}\label{eqquality}
f(x_N) - f(x_*) \leq \frac{R^2}{A_N} + \frac{2 \sum\limits_{k=0}^{N-1} \delta_{k+1}A_{k+1}}{A_{N}},
\end{equation}
где $A_N = \sum\limits_{k=0}^{N-1} \frac{1}{L_{k+1}}$.
\end{theorem}

Для получения оценки скорости сходимости алгоритма \ref{alg2} для задач оптимизации с относительной точностью выберем критерий остановки вида
\begin{gather}\label{krost}
A_N \geq \frac{R}{\varepsilon}.
\end{gather}
Покажем, как на базе предыдущей теоремы можно выбрать параметры алгоритма, чтобы получить аналог универсального градиентного метода для задач с относительной точностью вида \eqref{eq1}.

\begin{theorem}\label{then3} Пусть $f$~--- положительная однородная выпуклая функция, причём $f(x) \geq \gamma_0\|x\|_2$ для всякого $x \in Q$, а также $\|x_{*} - x_0 \|_2 \leq R\sqrt{2}\leq\|x_* - x_0\|_2 \sqrt{2}$ для некоторой постоянной $R > 0$, где $x_0:\, \|x_0\|_2 =\min\{\|x\|_2\;|\; x \in Q\}$~--- начальная точка, а $x_*$~--- ближайшая точка минимума $f$ к $x_0$ по евклидовой норме. Тогда после выполнения критерия остановки \eqref{krost} в случае $\varepsilon = \frac{\gamma \gamma_0}{3}$ и $\delta_{k+1} \leq R \varepsilon \frac{\alpha_{k+1}}{4 A_{k+1}}$ для всякого $k \geq 0$ будет гарантированно достигнута $\gamma$-относительная точность приближенного решения $x_{N}$: $f(x_{N}) \leq (1 + \gamma) f^*$. При этом критерий остановки \eqref{krost} будет выполнен не более, чем после
\begin{equation}\label{eq-complexity-alg2}
N \geq \left ( R \varepsilon^{-\frac{2}{1+\nu}} 2^{\frac{2 + 4 \nu}{1 + \nu}} L_{\nu}^{\frac{2}{1 + \nu}} \right )^{\frac{1+\nu}{1 + 3 \nu}}
\end{equation}
итераций алгоритма \ref{alg2}.
\end{theorem}
\begin{proof}
Выбор в алгоритме \ref{alg2} параметров
$$\delta_{k+1} \leq R \varepsilon \frac{\alpha_{k+1}}{ 4 A_{k+1}} \quad \forall k \geq 0$$
означает, что из теоремы~\ref{then2} вытекает неравенство
$$f(x_{N}) - f(x_{\ast}) \leq \frac{R^2}{A_N} + \frac{2 \sum_{k=0}^{N-1} A_{k+1} \delta_{k+1}}{A_N}.$$
Далее, имеет место неравенство (см. \cite{3}, п. 4.3, а также \cite{4})
$$A_N \geq \frac{N^{\frac{1+3 \nu}{1+ \nu}} \varepsilon^{\frac{1- \nu}{1+ \nu}}}{2^{\frac{2 + 4 \nu}{1 + \nu}} L_{\nu}^{\frac{2}{1+ \nu}}}.$$
Критерий остановки~\eqref{krost} алгоритма \ref{alg2} заведомо выполнен, если
$$\frac{N^{\frac{1+3 \nu}{1+ \nu}} \varepsilon^{\frac{1- \nu}{1+ \nu}}}{2^{\frac{2 + 4 \nu}{1 + \nu}} L_{\nu}^{\frac{2}{1+ \nu}}} \geq \frac{R}{\varepsilon}.$$
Следовательно,
$$N \geq \varepsilon^{\frac{-2}{1 + 3 \nu}} \left ( R \; 2^{\frac{2 + 4 \nu}{1 + \nu}} L_{\nu}^{\frac{2}{1+ \nu}} \right )^{\frac{1 + \nu}{1 + 3 \nu}}.$$

Пусть $\varepsilon = \frac{\gamma \gamma_0}{3}$. Тогда с учётом сделанного предположения о выборе параметров $\delta_{k+1}$ после выполнения критерия остановки \eqref{krost} алгоритма \ref{alg2} имеем неравенство:
$$f(x_N) - f^* \leq \frac{3 \varepsilon R}{2}.$$
Cогласно теореме~\ref{then1} верно $R \leq \frac{2}{\gamma_0} f^{\ast}$, откуда имеем неравенство
$$f(x_N) \leq f^* + \frac{3 \varepsilon}{\gamma_0} f^{*}.$$
Наконец, при $\varepsilon = \frac{\gamma \gamma_0}{3}$ получаем
$$f(x_N) \leq (1 + \gamma) f^*,$$
что и требовалось.
\end{proof}

Доказанная в предыдущей теореме оценка скорости сходимости алгоритма \ref{alg2}, как известно \cite{NemYudin}, оптимальна на выделенном классе задач минимизации выпуклых функционалов с гёльдеровым градиентом с точностью до умножения на постоянную величину. Однако за счёт универсальности метода (адаптивной настройки на параметр и уровень гладкости задачи $\nu \in [0; 1]$ в процессе работы) возможно ожидать на практике более высокую скорость повышения качества выдаваемого решения с ростом количества итераций. Для иллюстрации этого факта приведём пример с результатами некоторых вычислительных экспериментов.
\begin{example}\label{example1}
$f(x) = \|x\|_2, \; x \in \mathbb{R}^n, \; n = 10^6, \; Q = B_1 (2p)$ --- евклидов шар единичного радиуса с центром в точке $2p$, где $p = \frac{(1, \ldots, 1)}{\|(1,\ldots, 1)\|_2} \in \mathbb{R}^{10^6}$. Очевидно, что в таком случае $f^{*} > 0$ и $\gamma_0 = 1$. Отметим, что полученная в теореме \ref{then3} оценка скорости сходимости \eqref{eq-complexity-alg2} сублинейна. В лучшем случае ($\nu = 1$) оценка \eqref{eq-complexity-alg2} может гарантировать оценку сложности $O(\sqrt{\gamma^{-1}})$, которая неулучшаема на классе задач минимизации выпуклых функционалов с липшицевым градиентом. Однако, как видно по таблице 1 ниже, адаптивный подбор параметров $L_{k+1}$ ($k = 0, 1, ...$) приводит к повышению скорости сходимости.

\begin{table}[h]
\caption{Скорость изменения качества решения по мере роста числа итераций}
\begin{center}
\begin{tabular}{|c|c|c|}
\hline
Кол-во итераций & Теор. оценка качества реш. согл. \eqref{eqquality} & Время работы (сек.)\\
\hline
10 & 9.996e-05 &0:00:02.91\\
\hline
15 & 3.061e-06 & 0:00:04.46 \\
\hline
20 & 9.540e-08 & 0:00:05.91\\
\hline
25 & 2.980e-09 & 0:00:07.37\\
\hline
30 & 9.313e-11 & 0:00:08.84 \\
\hline
35 & 2.910e-12& 0:00:10.33 \\
\hline
40 & 9.094e-14 & 0:00:11.84 \\
\hline
45 & 2.842e-15 & 0:00:18.93 \\
\hline
50 & 8.881e-17 & 0:00:20.38 \\
\hline
55 & 2.775e-18 & 0:00:22.27 \\
\hline
60 & 8.673e-20 & 0:00:23.78 \\
\hline
65 & 2.710e-21 &0:00:25.31 \\
\hline
70 & 8.470e-23 & 0:00:26.90 \\
\hline
75 & 2.646e-24 & 0:00:28.39 \\
\hline
80 & 8.271e-26 & 0:00:29.98 \\
\hline
85 & 2.584e-27 & 0:00:31.49\\
\hline
90 & 8.077e-29& 0:00:32.99 \\
\hline
95 & 2.524e-30 & 0:00:34.47 \\
\hline
100 & 7.888e-32 & 0:00:35.90 \\
\hline
\end{tabular}
\end{center}
\end{table}

Вычисления были произведены с помощью Python 3.4 на компьютере с Intel(R) Core(TM) i7-8550U CPU @ 1.80GHz, 1992 Mhz, 4 Core(s). ОЗУ компьютера составляла 8 Гб.
\end{example}

Сделаем комментарий о совместимости предположения о гёльдеровости градиента \eqref{eqHolderGrad} целевой функции и естественным для задач с относительной точностью требованием положительной однородности $f$. Во-первых, если $0$ не лежит в допустимом множестве задачи $Q$, то на $Q$ целевая функция $f$ может иметь как гёльдеров, так и липшицев градиент, что верно к примеру для задачи из примера \ref{example1}. Более того, сама постановка задачи с относительной точностью подразумевает, что $f^{*} > 0$ и поэтому $0$ не должен лежать в допустимом множестве задачи. Во-вторых, липшицева функция $f$ удовлетворяет \eqref{eqHolderGrad} при $\nu = 0$, а адаптивность алгоритма \ref{alg2} потенциально позволяет улучшить оценку скорости сходимости по сравнению с оптимальной при $\nu = 0$ оценкой вида $O(\gamma^{-2})$, что как раз и показывают результаты экспериментов из примера \ref{example1}.

\section{Адаптивный субградиентный метод для минимизации слабо $\beta$-квазивыпуклых функционалов с острым минимумом}

Сублинейность оценки скорости сходимости в теореме~\ref{then3} приводит к идее исследовать подклассы задач, для которых всё же возможна линейная скорость сходимости субградиентного метода. Это, в частности, возможно в случае предположений о том, что минимум острый и доступно $f^{\ast}$ \cite{Pol_1969}.  Поэтому данный и следующий разделы работы посвящены исследованию методов с линейной скоростью сходимости для задач с острым минимумом на классах задач с обобщениями выпуклости (слабая $\beta$-квазивыпуклость и обычная квазивыпуклость). Отметим, что значение $f^{\ast}$ бывает известно в геометрических задачах (проекция точки на множество или нахождения общей точки системы множеств). Также можно упомянуть задачу типа $f(x)=\|Ax-b\|_{2}\rightarrow\min$, где $A$ --- матрица $n\times n$, $b\in\mathbb{R}^{n}$. Она разрешима, если существует $x_{*}$ такое, что $Ax_{*}=b$. Будем рассматривать задачи вида
\begin{equation}\label{eq_1}
f(x)\rightarrow\min_{x\in Q},
\end{equation}
где $Q$~--- выпуклое замкнутое подмножество $\mathbb{R}^{n}$, $f^{*} = f(x_*) = \min\limits_{x \in Q} f(x)$, а $f$~--- слабо $\beta$-квазивыпуклый функционал при некотором $\beta\in (0;1]$. Напомним (\cite{Hardt}), что $f$ называется {\it слабо $\beta$-квазивыпуклым} относительно точки минимума $x_{*}$ задачи \eqref{eq_1} на множестве $Q$, если для произвольного $x\in Q$ выполнено неравенство:
\begin{equation}\label{eqquasiconv}
f(x_{*})\geqslant f(x)+\frac{1}{\beta} \langle \nabla f(x), x_{*}-x \rangle,
\end{equation}
где $\nabla f(x)$~--- произвольный субградиент $f$ в точке $x$. Под субградиентом мы здесь и всюду далее понимаем элемент субдифференциала Кларка $f$ в точке $x$ и предполагаем его существование. Если $f$ дифференцируем в точке $x$, то под $\nabla f$ понимаем обычный градиент. Это вполне естественно для липшицевых функционалов (для существования субдифференциала Кларка в точке достаточно локальной липшицевости $f$ в окрестности этой точки). Если функционал $f$ выпуклый, то субдифференциал Кларка совпадает с обычным субдифференциалом в смысле выпуклого анализа, и в таком случае условие слабой $\beta$-квазивыпуклости \eqref{eqquasiconv} верно при $\beta = 1$.

Ясно, что если неравенство \eqref{eqquasiconv} верно для некоторого $\beta = \beta_0 \in (0; 1]$, то оно верно и при $\beta \in (0; \beta_0]$. Примеры функционалов, для которых возможно проверить свойство слабой $\beta$-квазивыпуклости и оценить параметр $\beta$, приведены в \cite{Hind}. В частности, это верно для невыпуклой функции $f(x)=|x|(1-e^{-|x|})$ при $\beta=1$ \cite{Hardt}. Отметим также, что субградиент $f$ может быть нулевым только в точке минимума: равенство $\nabla f(x)=0$ влечет $f(x)\leqslant f(x_*)$, что автоматически означает $f(x)=f^{*}$.

Предположим, что также верно условие острого минимума для некоторого $\alpha >0$
\begin{equation}\label{tochka}
f(x) - f^* \geqslant \alpha\cdot \min_{x_* \in X_*} \|x - x_*\|_{2},
\end{equation}
где $X_*$~--- множество точек минимума функции $f$ на множестве $Q$. В частности, условие \eqref{tochka} верно для задачи евклидова проектирования точки $x$ на выпуклый компакт $X_{*}\subset Q$, причём $f^{*}=0$.

Далее будем понимать под $Pr_Q$ оператор евклидова проектирования на множество $Q$. Предложим следующий вариант субградиентного метода с шагом Б.\,Т.\,Поляка \cite{Pol_1969} для задачи минимизации слабо $\beta$-квазивыпуклого функционала:
\begin{equation}\label{1}
x_{k+1} = Pr_{Q} \{x_k - h_k \nabla f(x_k)\},
\end{equation}
где при всяком $k \geqslant 0$ верно $\|\nabla f(x_k) \|_2 \neq 0$ (иначе решение уже найдено), а также
$$h_k = \frac{\beta(f(x_k) - f^{\ast})}{\| \nabla f(x_k) \|_2^2}.$$

Справедлива следующая

\begin{theorem}\label{thmbetaquasiconv}
Пусть $f$~--- слабо $\beta$-квазивыпуклая функция и для задачи минимизации $f$ \eqref{eq_1} с острым минимумом используется алгоритм \eqref{1} c шагом
$$h_k = \dfrac{\beta(f(x_k) - f^{\ast})}{\| \nabla f(x_k) \|_2^2}.$$
Тогда после $k$ итераций алгоритма \eqref{1} верно неравенство:
$$\min_{x_* \in X_*}\|x_{k+1} - x_* \|_2^2 \leq \prod_{i=0}^k \left ( 1 - \frac{\alpha^2\beta^2}{\| \nabla f(x_i)\|_2^2} \right) \min_{x_* \in X_*}\|x_0 - x_* \|_2^2.$$
\end{theorem}
\begin{proof}
Как известно (\cite{Pol_1969}, соотношения (3) из доказательства теоремы 1), для всякого $x_* \in X_*$ верны неравенства
$$2\beta h_k (f(x_k) - f(x_*)) \leq 2h_k \langle \nabla f(x_k), x_k-x_* \rangle \leq h_k^2 \| \nabla f(x_k)\|_2^2 + \min_{x_* \in X_*}\|x_k - x_* \|_2^2 - \min_{x_* \in X_*}\|x_{k+1} - x_* \|_2^2 .$$
Поэтому
$$\min_{x_* \in X_*} \|x_{k+1} - x_* \|_2^2 \leq h_k^2 \| \nabla f(x_k)\|_2^2 - 2 \beta^{2} h_k (f(x_k) - f(x_*)) + \min_{x_* \in X_*} \|x_k - x_* \|_2^2 = $$
$$ = \frac{\beta^2(f(x_k) - f(x_*))^2}{\| \nabla f(x_k) \|_2^2} - \frac{2\beta^2(f(x_k) - f(x_*))^2}{\| \nabla f(x_k) \|_2^2} + \min_{x_* \in X_*} \|x_k - x_* \|_2^2 =
$$
$$= - \frac{\beta^2(f(x_k) - f(x_*))^2}{\| \nabla f(x_k) \|_2^2} + \min_{x_* \in X_*} \|x_k - x_* \|_2^2 .$$
Согласно условию острого минимума имеем следующие соотношения:
$$\min_{x_* \in X_*} \|x_{k+1} - x_* \|_2^2 \leq - \frac{\alpha^2\beta^2}{\| \nabla f(x_k)\|_2^2} \min_{x_* \in X_*} \|x_k - x_* \|_2^2 + \min_{x_* \in X_*}\|x_k - x_* \|_2^2 =
$$
$$
= \left ( 1 - \frac{\alpha^2\beta^2}{\| \nabla f(x_k)\|_2^2} \right) \min_{x_* \in X_*}\|x_k - x_* \|_2^2.
$$
Далее, получаем цепочку неравенств:
$$\min_{x_* \in X_*} \|x_{k+1} - x_* \|_2^2 \leq \left ( 1 - \frac{\alpha^2\beta^2}{\| \nabla f(x_k)\|_2^2} \right) \min_{x_* \in X_*}\|x_k - x_* \|_2^2 \leq $$
$$ \leq \left ( 1 - \frac{\alpha^2\beta^2}{\| \nabla f(x_k)\|_2^2} \right) \left ( 1 - \frac{\alpha^2\beta^2}{\| \nabla f(x_{k-1})\|_2^2} \right) \min_{x_* \in X_*} \|x_k - x_* \|_2^2 \leq $$
$$ \leq \ldots \leq \prod_{i=0}^k \left ( 1 - \frac{\alpha^2\beta^2}{\| \nabla f(x_i)\|_2^2} \right) \min_{x_* \in X_*}\|x_0 - x_* \|_2^2.$$
\end{proof}

\begin{corollary}\label{corol}
Если в условиях предыдущей теоремы допустить, что $f$ удовлетворяет условию Липшица с константой $M >0$, то можно утверждать сходимость алгоритма \eqref{1} со скоростью геометрической прогрессии:
$$\min_{x_* \in X_*}\|x_{k+1} - x_* \|_2^2 \leq \left ( 1 - \frac{\alpha^2\beta^2}{M^2} \right)^{k+1} \min_{x_* \in X_*}\|x_0 - x_* \|_2^2.$$
\end{corollary}

Отметим, что теорема \ref{thmbetaquasiconv} применима к функционалам, не удовлетворяющим условию Липшица. Это возможно ввиду адаптивности оценки скорости сходимости из теоремы \ref{thmbetaquasiconv} (использование норм субградиентов $\| \nabla f(x_i)\|_2^2$ вместо $M^2$). Но если нет уверенности в выполнении условия Липшица, то в данной общей ситуации невозможно гарантировать линейную скорость сходимости.

Покажем, как можно получить оценку скорости сходимости алгоритма \eqref{1} для достижения $\gamma$-относительной точности по целевому функционалу
$f(x_{k+1})\leq (1 + \gamma) f^*$ при $\gamma>0$. Следуя работам~\cite{Nest_diss, NesRound, NesUnco}, будем рассматривать задачу с относительной точностью на классе выпуклых ($\beta = 1$) положительно однородных функций $f:\, Q\rightarrow \mathbb{R}:\,f^{*}>0$, а также с условием $f(x) \geq \gamma_0 \|x\|_2$ для всякого $ x \in Q$ при некотором фиксированном $\gamma_0 >0$. Начальная точка $x_0$ рассматриваемого алгоритма \eqref{1} выбирается так, чтобы $\| x_0 \|_2 := \min_{x \in Q}\{\|x\|_2\}$.

По следствию \ref{corol} имеем неравенства:
$$\min_{x_{*}\in X_{*}}\|x_{k+1} - x_*\| \leq \left ( 1 -\frac{\alpha^2}{M^2} \right )^{k+1}\min_{x_{*}\in X_{*}} \|x_0 - x_* \|_2^2 .$$
С учётом теоремы \ref{then1} имеем:
$$\|x_{k+1} - x_* \| \leq \frac{2 f^*}{\gamma_0} \sqrt{\left ( 1 -\frac{\alpha^2}{M^2} \right )^{k+1}} \leq \gamma f^*,$$
$$f(x_{k+1}) - f^* \leq \frac{2 M f^*}{\gamma_0} \sqrt{\left ( 1 -\frac{\alpha^2}{M^2} \right )^{k+1}} \leq \gamma f^*,$$
$$f(x_{k+1}) \leq f^* + \gamma f^* = (1 + \gamma) f^* .$$
Потребуем, чтобы
$$2 M \sqrt{\left ( 1 -\frac{\alpha^2}{M^2} \right )^{k+1}} \leq \gamma \gamma_0 .$$
Это означает, что $$k \geq \log_{\left ( 1 -\frac{\alpha^2}{M^2} \right )} \frac{\gamma^2 \gamma_0^2}{4 M^2} - 1 = \frac{2 (\ln{2M} - \ln{\gamma \gamma_0})}{\ln{M^2} - \ln{(M^2 - \alpha^2)}} - 1.$$

Таким образом, из следствия \ref{corol} при $\beta = 1$ (в выпуклом случае) имеем
\begin{corollary}\label{corol2}
Пусть $f$~--- выпуклая $M$-липшицева положительно однородная функция, причём для всякого $x\in Q$ верно неравенство $\gamma_0\|x\|_2 \leq f(x)$. Тогда для всякого достаточно малого $\gamma>0$ при условии
$\|x_0\|_2 = \min_{x \in Q} \|x\|_2$ после
$$k \geq \frac{2 (\ln{2M} - \ln{\gamma \gamma_0})}{\ln{M^2} - \ln{(M^2 - \alpha^2)}} - 1$$
итераций алгоритма \eqref{1} имеем неравенство
$$f(x_{k+1})\leq (1 + \gamma) f^{*}.$$
\end{corollary}

\section{Об универсальности одного субградиентного метода для задач минимизации гёльдеровых квазивыпуклых функционалов с острым минимумом}

Теперь рассмотрим класс задач \eqref{eq_1} минимизации квазивыпуклых функционалов $f$
$$
f(\lambda x+(1-\lambda)y)\leq\max\{f(x),f(y)\}\quad \forall\,x,y\in Q,\, \lambda\in[0;1].
$$
в предположении о том, что $f$ удовлетворяет условию Гёльдера
\begin{equation}\label{eeq_2}
|f(x)-f(y)|\leq M_{\nu}\|x-y\|^{\nu}_{2}\quad\forall\,x,y\in Q
\end{equation}
при некотором фиксированном $\nu\in[0;1]$, $0 \leq M_{\nu} < + \infty$.

Введём (следуя~\cite{NesEffec}, \cite{Nest2010mon}) вспомогательную величину
$$
\upsilon_{f}(x, x_{*}) :=
\begin{cases}
\left\langle\frac{\nabla f(x)}{\|\nabla f(x)\|_{2}},x-x_{*}\right\rangle, & \text{если $x\neq x_{*}$;} \\
0 & \text{при $x= x_{*}$}
\end{cases}
$$
Если $\nabla f(x)=0$ при $x\neq x_{*}$, то вместо $\nabla f(x)$ можно использовать ненулевой вектор нормали $\mathcal{D}f(x)$ ко множеству уровня функции $f$ в точке $x$~\cite{NesEffec}. Но для упрощения изложения далее сделаем допущение, что $\nabla f(x)\neq0$ при $x\neq x_{*}$.

С учётом \eqref{tochka} и \eqref{eeq_2} получаем, что для всякого $x\in Q$
$$
\alpha\min_{x_{*}\in X_{*}}\|x-x_{*}\|_{2}\leq f(x)-f^{*}\leq M_{\nu}\min_{x_{*}\in X_{*}}\|x-x_{*}\|_{2}^{\nu}.
$$
Поэтому
$$
\min_{x_{*}\in X_{*}}\|x-x_{*}\|_{2}^{1-\nu}\leq\frac{M_{\nu}}{\alpha},
$$
откуда при $\nu<1$
$$
\min_{x_{*}\in X_{*}}\|x-x_{*}\|_{2}\leq \left(\frac{M_{\nu}}{\alpha}\right)^{\frac{1}{1-\nu}}.
$$
Поэтому при $0\leq\nu<1$ можно локализовать допустимую область $Q$, заменив её на пересечение с евклидовым шаром с центром в точке $x_{0}$ и радиусом $\displaystyle\left(\frac{M_{\nu}}{\alpha}\right)^{\frac{1}{1-\nu}}$.

Справедлива следующее

\begin{propose}\label{Lemma1}
Если $f$ квазивыпукла и удовлетворяет \eqref{eeq_2}, то существует такое $M>0$, что для всякого $x\in Q$
\begin{equation}\label{eeq_3}
f(x)-f(x_{*})\leq M\upsilon_{f}(x,x_{*}),
\end{equation}
\end{propose}
\begin{proof}
Случай $\upsilon_{f}(x,x_{*})=0$ очевиден.
Неравенство \eqref{eeq_2} означает, что
\begin{equation}\label{eeq_4}
f(x)-f(x_{*})\leq M_{\nu}\upsilon_{f}^{\nu}(x,x_{*})\quad\forall\,x\in Q
\end{equation}
ввиду леммы 3.2.1 из монографии ~\cite{Nest2010mon}. Если в \eqref{eeq_4} $\upsilon_{f}(x,x_{*})\geq 1$ для некоторого $x\in Q$, то
$$
f(x)-f(x_{*})\leq M_{\nu}\upsilon_{f}(x,x_{*}).
$$

Если же  $\nu<1$ и $\upsilon_{f}(x,x_{*})< 1$, то применим компактность допустимого множества $Q$. Субдифференцируемость $f$ в произвольной точке $Q$ означает локальную липшицевость $f$. Ввиду компактности $Q$ это означает $M_{f}$--липшицевость $f$ для некоторой $M_{f}>0$. Остаётся теперь обратить внимание на \eqref{eeq_4} при $\nu=1$.
\end{proof}

\begin{remark}
Если $a$ не является сильно малым, то верно неравенство (аналогичное неравенству из замечания 5.1 из \cite{Gasnikov_book})
$$
M_{\nu}a^{\nu}\leq\frac{M_{\nu}^{\frac{2}{1+\nu}}}{2}\frac{a^{2}}{\delta^{\frac{1-\nu}{1+\nu}}}+\frac{\delta}{2}\qquad\forall\,\delta>0.
$$
Такое неравенство позволит оценить
$$
M=M(0)=\max\left\{M_{\nu},\frac{M_{\nu}^{\frac{2}{1+\nu}}}{2}\right\}
$$
при достаточно большом $\upsilon_{f}(x,x_{*})$ для некоторого $x_{*}\in X_{*}$.
\end{remark}

Рассмотрим метод $(k=0,1,2,\ldots)$
\begin{equation}\label{eeq_5}
x_{k+1}=Pr_{Q}\{x_{k}-h_{k}\nabla f(x_{k})\},\,\text{где}\,h_{k}=\frac{f(x_{k})-f(x_{*})}{M\|\nabla f(x_{k})\|_{2}},
\end{equation}
где $M$ удовлетворяет \eqref{eeq_3}. Докажем, что данный алгоритм универсален в том смысле, что сходится со скоростью геометрической прогрессии для задачи минимизации квазивыпуклого гёльдерова функционала с острым минимумом и заранее известным точным значением $f^{*}$ при всех $\nu \in [0; 1]$.

\begin{theorem}\label{th5}
Пусть верно \eqref{eeq_3} при некотором $0 < M < + \infty$ и $f$ имеет $\alpha$-острый минимум \eqref{tochka}. Тогда метод \eqref{eeq_5} сходится со скоростью геометрической прогрессии:
\begin{equation}\label{eeq_6}
\min_{x_* \in X_*}\|x_{k+1}-x_{*}\|^{2}_{2}\leq\left(1-\frac{\alpha^{2}}{M^{2}}\right)^{k+1} \min_{x_* \in X_*}\|x_{0}-x_{*}\|^{2}_{2}.
\end{equation}
\end{theorem}
\begin{proof}
Действительно, для метода \eqref{eeq_5} верны неравенства (см. \cite{Pol_1969}, соотношения (3) из доказательства теоремы 1)
$$
2h_{k}\langle\nabla f(x_{k}),x_{k}-x_{*}\rangle\leq h_{k}^{2}\|\nabla f(x_{k})\|^{2}_{2}+\|x_{k}-x_{*}\|^{2}_{2}-\|x_{k+1}-x_{*}\|^{2}_{2}
\quad\forall\, k\geq0.
$$
Поэтому имеем:
$$
\min_{x_* \in X_*}\|x_{k+1}-x_{*}\|^{2}_{2}\leq\min_{x_* \in X_*}\|x_{k}-x_{*}\|^{2}_{2}+ h^{2}_{k}\|\nabla f(x_{k})\|^{2}_{2} - 2h_{k}\langle\nabla f(x_{k}),x_{k}-x_{*}\rangle=
$$
$$
=\min_{x_* \in X_*}\|x_{k}-x_{*}\|^{2}_{2}+\frac{(f(x_{k})-f(x_{*}))^{2}}{M^{2}}-
\frac{2(f(x_{k})-f(x_{*}))}{M}\cdot\left\langle\frac{\nabla f(x_{k})}{\|\nabla f(x_{k})\|_{2}},x_{k}-x_{*}\right\rangle=
$$
$$
=\min_{x_* \in X_*}\|x_{k}-x_{*}\|^{2}_{2}+ \left(\frac{f(x_{k})-f(x_{*})}{M}\right)^{2}-\frac{2(f(x_{k})-f(x_{*}))}{M} \cdot
\upsilon_{f}(x_{k}, x_{*}),
$$
$$
\leq \min_{x_* \in X_*}\|x_{k}-x_{*}\|^{2}_{2} - \left(\frac{f(x_{k})-f(x_{*})}{M}\right)^{2}\leq
\min_{x_* \in X_*}\|x_{k}-x_{*}\|^{2}_{2}
\left(1-\frac{\alpha^{2}}{M^{2}}\right)
$$
ввиду \eqref{tochka}, т.к. $\upsilon_{f}(x_{k}, x_{*})\geq\frac{f(x_{k})-f(x_{*})}{M}$. Ясно, что
$$
\min_{x_* \in X_*}\|x_{k+1}-x_{*}\|^{2}_{2}\leq
\left(1-\frac{\alpha^{2}}{M^{2}}\right)\min_{x_* \in X_*}\|x_{k}-x_{*}\|^{2}_{2}\leq
$$
$$
\leq \left(1-\frac{\alpha^{2}}{M^{2}}\right)\min_{x_* \in X_*}\|x_{k-1}-x_{*}\|^{2}_{2}\leq\ldots\leq
\left(1-\frac{\alpha^{2}}{M^{2}}\right)^{k+1}\min_{x_* \in X_*}\|x_{0}-x_{*}\|^{2}_{2},
$$
что и требовалось.
\end{proof}

Сформулируем теперь следствие из теоремы \ref{th5} об оценке скорости сходимости алгоритма \eqref{eeq_5} для задач выпуклой положительно однородной минимизации с относительной точностью по целевом функционалу.

\begin{corollary}\label{corol3}
Пусть $f$~--- выпуклая $M$-липшицева положительно однородная функция, причём для всякого $x\in Q$ верно неравенство $\gamma_0\|x\|_2 \leq f(x)$. Тогда для всякого достаточно малого $\gamma>0$ при условии
$\|x_0\|_2 = \min_{x \in Q} \|x\|_2$ после
$$k \geq \frac{2 (\ln{2M} - \ln{\gamma \gamma_0})}{\ln{M^2} - \ln{(M^2 - \alpha^2)}} - 1$$
итераций алгоритма \eqref{eeq_5} имеем неравенство $f(x_{k+1})\leq (1 + \gamma) f^{*}.$
\end{corollary}

\begin{remark}\label{remdeltasharmpmin}
Некоторым недостатком алгоритма \eqref{eeq_5} по сравнению с алгоритмом \eqref{1} может считаться требование знать $M$ из неравенства леммы \ref{Lemma1} при организации шагов. Однако эта же особенность позволяет предложить вариацию субградиентного метода для задач с некоторым $\Delta$-обобщением острого минимума вида
$$
f(x)-\overline{f}\geq\alpha\min_{x_{*}\in X_{*}}\|x-x_{*}\|_{2}-\Delta\quad\forall\,x\in Q
$$
при фиксированном значении $\Delta>0$ и заданном $\overline{f}\geq f^{*}$. В частности, данное условие логично использовать при отсутствии точной информации об $f^{*}$. Если выбрать в алгоритме \eqref{eeq_5} шаги вида
$$
h_{k}=\frac{f(x_{k})-\overline{f}}{M\|\nabla f(x_{k})\|_{2}},
$$
то при $2M^2 > \alpha^2$ получим соотношения (аналогично рассуждениям из доказательства теоремы \ref{th5})
$$
\min_{x_{*}\in X_{*}}\|x_{k+1}-x_{*}\|_{2}^{2}\leq \min_{x_{*}\in X_{*}}\|x_{k}-x_{*}\|_{2}^{2}-\frac{1}{M^{2}}(f(x_{k})-\overline{f})^{2}\leq
$$
$$
\left[\text{применим неравенство}\quad a^{2}\geq \frac{1}{2}(a+b)^{2}-b^{2},\text{ справедливое для всяких }\,a,b\in\mathbb{R}\right]
$$
$$
\leq \min_{x_{*}\in X_{*}}\|x_{k}-x_{*}\|_{2}^{2}\left(1-\frac{\alpha^{2}}{2M^{2}}\right)+\frac{\Delta^{2}}{M^{2}}\leq
$$
$$
\leq \min_{x_{*}\in X_{*}}\|x_{0}-x_{*}\|_{2}^{2}\left(1-\frac{\alpha^{2}}{2M^{2}}\right)^{k+1}+\frac{\Delta^{2}}{M^{2}}\left(1+\left(1-\frac{\alpha^{2}}{2M^{2}}\right)+\ldots+\left(1-\frac{\alpha^{2}}{2M^{2}}\right)^{k}\right)\leq
$$
$$
\leq \min_{x_{*}\in X_{*}}\|x_{0}-x_{*}\|_{2}^{2}\left(1-\frac{\alpha^{2}}{2M^{2}}\right)^{k+1}+\frac{2\Delta^{2}}{\alpha^{2}}.
$$
Таким образом, наличие неточности $\Delta$ приводит к дополнительному слагаемому вида $O(\Delta^{2})$ в оценке невязки $\min_{x_* \in X_*} \|x_{k+1} - x_*\|_2$. По-видимому, в случае малых значений $\|\nabla f(x_k)\|_2$ для алгоритма \eqref{1} такой вывод сделать уже нельзя.
\end{remark}

На первый взгляд условие \eqref{tochka} представляется довольно ограничительным. Однако это не совсем так, что показано в следующих примерах.

\begin{example}
Пусть $f$ имеет {\sl квадратичный рост}:
\begin{equation}\label{ex1}
f(x)-f(x_{*})\geq\mu\|x-x_{*}\|_{2}^{2}\quad\forall\,x\in Q
\end{equation}
при некотором $\mu > 0$. Например, условие \eqref{ex1} заведомо верно для всякой $2\mu$-сильно выпуклой $f$.
Условие \eqref{ex1} означает, что
$$
F(x):=\sqrt{f(x)-f(x_{*})}\geq\sqrt{\mu}\|x-x_{*}\|_{2}
$$
при всяком $x\in Q$ и $F$ удовлетворяет условию острого минимума \eqref{tochka} при $\alpha=\sqrt{\mu}$. Если $f$ квазивыпукла, то и $F$ квазивыпукла. Если $f$ $M_{f}$--липшицева, то неравенство $|\sqrt{a}-\sqrt{b}|\leq\sqrt{|a-b|}\quad(a,b\in\mathbb{R})$ влечёт гёльдеровость $F$ при $\nu=\dfrac{1}{2}$ (аналогичный вывод можно сделать и для гёльдеровых $f$). Это приводит к выводу: оценка \eqref{eeq_6} верна для минимизации $M_{f}$--липшицевой квазивыпуклой функции $f$ с квадратичным ростом в предположении $f(x_{*})$.
\end{example}

\begin{example}\label{example3}
Во многих прикладных задачах возникает также условие так называемого {\sl слабого острого минимума} (в западной литературе для этого понятия можно встретить термин {\sl условие гёльдерова роста}; см., например, \cite{Jiang, Johnstone, LiuYang})
$$
f(x)-f^{*}\geq\mu\min_{x_* \in X_*}\|x-x_{*}\|_{2}^{p}\quad\forall\,x\in Q
$$
для некоторого фиксированного $p\in [1; + \infty)$ и некоторой постоянной $\mu>0$. Ясно, что это условие обобщает как обычное условие острого минимума~\eqref{tochka}, так и условие квадратичного роста~\eqref{ex1}. В этом случае снова возможно применить схему рассуждений предыдущего примера.
Действительно, если $f$ квазивыпукла, то функция
$$
F(x):=(f(x)-f^{*})^{1/p}
$$
также квазивыпукла и $\min\limits_{x\in Q}F(x)=0$, т.е. $F$ имеет $(\mu^{1/p})$--острый минимум. В случае, если $f$ удовлетворяет условию Гёльдера с показателем $\nu\in[0;1]$, то $F$ также удовлетворяет условию Гёльдера с показателем $\dfrac{\nu}{p}$, причём $\dfrac{\nu}{p}\in[0;1]$. Таким образом, к задаче минимизации $f$ применима оценка~\eqref{eeq_6} теоремы~\ref{th5}, означающая линейную скорость сходимости. Это указывает на некоторую универсальность алгоритма~\eqref{eeq_5} по параметру $p \in [1; + \infty)$. Отметим, что ранее для задач со слабым острым минимумом при $\nu \neq p > 1$ были известны результаты только о сублинейной скорости сходимости субградиентных методов и вообще не рассматривался случай $\nu = 0$ (см. \cite{Rates} и имеющиеся в этой работе ссылки).
\end{example}

\section{Заключение}

В настоящей статье описаны новые результаты о скорости сходимости методов градиентного типа для двух вариантов постановки задачи. Во-первых, получена оценка скорости сходимости адаптивных методов для задач выпуклой однородной минимизации с гарантией достижения заданной относительной точности по целевому функционалу. Вторая часть статьи (разделы 3 и 4) посвящена результатам о линейной скорости сходимости субградиентных методов для задач с острым минимумом. Наиболее тонкий, на наш взгляд, результат получен в разделе 4 для класса задач минимизации квазивыпуклых гёльдеровых функционалов с острым минимумом. Как следствие, обоснована возможность построения субградиентного метода с линейной скоростью сходимости по аргументу для задач со слабым острым минимумом (гёльдеровым ростом) в предположении доступности информации о минимальном значении $f^{\ast}$. Последнее условие представляется довольно ограничительным, особенно для выводов по задачам с относительной точностью (следствия \ref{corol2} и \ref{corol3}). Поэтому в качестве возможного развития настоящей работы представляется интересным исследование вариаций предложенных методов с неточной информацией как об $f^{\ast}$ (см. замечание \ref{remdeltasharmpmin}), так и о значениях функции $f$, субградиентов $\nabla f$ в произвольных запрашиваемых точках допустимого множества решаемой задачи, то есть для задач с какими-то подходящими аналогами $(\delta,L)$-оракула~\cite{Devolder}. Авторы благодарят А. В. Гасникова и Ю. Е. Нестерова за полезные обсуждения.


\end{document}